\documentclass[review]{elsarticle}

\usepackage{lineno,hyperref}

\usepackage{epsfig}
\usepackage{amsmath}
\usepackage{amsfonts}
\usepackage{amsthm}
\usepackage{amsbsy}
\usepackage{appendix}
\usepackage{epstopdf}
\usepackage[table]{xcolor}
\usepackage{subfloat}
\usepackage{float}
\usepackage[thinlines]{easytable}
\usepackage{array}
\usepackage{color}
\usepackage{inputenc}
\usepackage{hyperref}

\numberwithin{equation}{section}
\allowdisplaybreaks

\theoremstyle{plain}

\newtheorem{theo}{Theorem}[section]
\newtheorem{algo}{Algorithm}[section]
\newtheorem{assump}{\bf Assumption}[section]

\newtheorem{definition}{\bf Definition}[section]
\theoremstyle{remark}

\theoremstyle{remark}

\modulolinenumbers[5]

\journal{Journal of \LaTeX\ Templates}

\begin{document}

\begin{frontmatter}

\title{Lagrangian Dynamic Mode Decomposition for Construction of Reduced-Order Models of Advection-Dominated Phenomena }
\author[mymainaddress]{Hannah Lu\fnref{email1}}
\fntext[email1]{email: hannahlu{@}stanford.edu}
\author[mymainaddress]{Daniel M. Tartakovsky\corref{mycorrespondingauthor}}
\cortext[mycorrespondingauthor]{Corresponding author}
\ead{tartakovsky@stanford.edu}
\address[mymainaddress]{Department of Energy Resources Engineering, Stanford, CA 94305. USA}

\begin{abstract}
Proper orthogonal decomposition (POD) and dynamic mode decomposition (DMD) are two complementary singular-value decomposition (SVD) techniques that are widely used to construct reduced-order models (ROMs) in a variety of fields of science and engineering. Despite their popularity, both DMD and POD struggle  to formulate accurate ROMs for advection-dominated problems because of the nature of SVD-based methods. We investigate this shortcoming of conventional POD and DMD methods formulated within the Eulerian framework. Then we propose a Lagrangian-based DMD method to overcome this so-called translational issues. Our approach is consistent with the spirit of physics-aware DMD since it accounts for the evolution of characteristic lines. Several numerical tests are presented to demonstrate the accuracy and efficiency of the proposed Lagrangian DMD method. 
\end{abstract}


\end{frontmatter}

\linenumbers
\section{Introduction}
Advection-diffusion equations are routinely used as a high-fidelity representation of mass conservation at a variety of spatiotemporal scales in a plethora of applications~\cite{tartakovsky-2019-diffusion}. These equations become highly nonlinear when advection velocity and/or diffusion coefficient depend(s) on a system state, e.g., in the case of multiphase flows in porous media~\cite{tartakovsky-2019-diffusion}. High-dimensional complex dynamics described by such nonlinear advection-diffusion equations often posses low-dimensional structures, suggesting the possibility of their replacement with reduced-order models (ROMs)~\cite{benner2015survey, quarteroni2014reduced, %ganapathysubramanian2007modeling,
acharjee2006concurrent%,ma2010adaptive
}.

Singular-value decomposition (SVD) can be utilized to extract a low-dimensional structure from the data generated with a high-fidelity model (HFM), i.e., to construct a conventional ROM. An example of such a ROM is built by deploying Galerkin projection to map a HFM onto a much smaller subspace; projection-based ROMs are referred to as proper orthogonal decomposition (POD)~\cite{holmes2012turbulence, lumley2007stochastic,  volkwein2011model}. The efficiency and accuracy of POD in nonlinear setting are increased by combining it with either the empirical interpolation method (EIM)~\cite{barrault2004empirical} or discrete empirical interpolation method (DEIM)~\cite{chaturantabut2010nonlinear}. Another example is the  Dynamic Mode Decomposition (DMD) method~\cite{kutz2016dynamic, mezic2013analysis}, which is used to discover a spatiotemporal coherent structure in the HMF-generated data. DMD's connection to the Koopman operator theory of nonlinear dynamic systems~\cite{koopman1931hamiltonian} is of theoretical interest~\cite{mezic2005spectral, mezic2013analysis, mezic2004comparison, rowley2009spectral}, while its equation-free spirit facilitates its use in conjunction with machine learning techniques~\cite{brunton2016discovering,schmidt2009distilling}.

While the robustness of DMD for parabolic problems has been established (including numerical~\cite{duke2012error} and theoretical~~\cite{korda2018convergence, lu2019error} analysis of its accuracy and convergence), both DMD and POD are known to fail in translational problems, such as wave-like phenomena, moving interfaces and moving shocks~\cite{kutz2016dynamic}.It can be explained by the intuition that the dominating advection behavior is traveling  through the whole high-dimensional domain, making it  impossible to determine a global spatiotemporal basis confined in a low-dimensional subspace. We facilitate this intuitive explaination with a concrete example in section~\ref{sec:Euler}. In terms of the Koopman operator theory, important physical observables (e.g., advection speed, shock speed, shock formation time) are unaccounted for in the standard DMD algorithm. Remedies for POD include the deployment of local basis~\cite{amsallem2012nonlinear}, domain decomposition~\cite{lucia2001reduced}, or basis splitting~\cite{carlberg2015adaptive}. A similar extension of DMD consists of multi-resolution DMD~\cite{kutz2016multiresolution}, which separates frequencies of different scales by filtering windows. Unfortunately, these remedies often compromise the ROM's efficiency by increasing its computational complexity. Alternative generalizations of DMD and POD explore symmetry and self-similarity properties to eliminate the translational issue using analytical tools~\cite{gerbeau2014approximated, kavousanakis2007projective, lumley2007stochastic, rapun2010reduced, rowley2003reduction, rowley2000reconstruction}. However, such tools are usually problem-dependent and mostly applicable to single-wave dominated problems.

Motivated by the recent work on Lagrangian POD~\cite{mojgani2017lagrangian}, we propose a physics-aware DMD method to construct a ROM within the Lagrangian framework. We chose the temporally evolving characteristic lines, a crucial physical quantity in advection-dominated systems, as a key observable of the underlying Koopman operator. Then, the DMD algorithm is used to identify, from sufficient data, a low-dimensional structure in the Lagrangian framework and thus to construct a physics-based ROM by approximating the underlying Koopman operator. The Lagrangian DMD can be applied to general advection-diffusion nonlinear flows. Furthermore, DMD outperforms POD in terms of computational costs due to the feature of iteration free. With the error analysis in \cite{lu2019error}, one can also estimate the accuracy of the ROM.

In Section~\ref{sec:Euler}, we review conventional POD and DMD in the Eulerian framework and demonstrate the issue of translational problems using advection-dominated linear flow as an example. Section~\ref{sec:Lagrange} provides a brief illustration of the Lagrangian POD in~\cite{mojgani2017lagrangian}, introduces our new Lagrangian DMD, and demonstrates its connection to the Koopman operator. Section~\ref{sec:example} contains several computational experiments used to validate the accuracy and efficiency of the proposed approach to ROM construction. It also compares the Lagrangian POD and DMD in terms of their accuracy and computational costs. Finally, in Section~\ref{sec:concl}, we draw conclusions and discuss the related ongoing work.

%%%%%%%%%%%%%%%%%%%%%%%%%%%%%%%%%%%%%%%
\section{Conventional Eulerian Reduced-Order Models}
\label{sec:Euler}
%%%%%%%%%%%%%%%%%%%%%%%%%%%%%%%%%%%%%%%

Consider a scalar state variable $u(x,t): [a,b] \times [0,T] \rightarrow \mathbb R^+$, whose dynamics is described by a one-dimensional nonlinear advection-diffusion equation
\begin{align}\label{2-1}
\frac{\partial u}{\partial t}+f(u)\frac{\partial u}{\partial x} = \frac{\partial}{\partial x}\left(D(x,t,u)\frac{\partial u}{\partial x}\right),\qquad f(u)=\frac{\partial F(u)}{\partial u},
\end{align}
subject to the initial condition
\begin{align}
u(x,t=0) = u_0(x)
\end{align}
and appropriate (arbitrary) boundary conditions at $x = a$ and $x = b$.
Within the Eulerian framework, the space is fixed and the interval $[a,b]$ is discretized with a uniform grid $\bold x = [x_1 = a,x_2\cdots,x_{N-1}, x_N = b]^\top$ of mesh size $\Delta x \equiv x_{j+1}-x_j = (b-a)/N$ and $N$ nodes. Likewise, the time interval $[0,T]$ is discretized uniformly with time step $\Delta t \equiv t^{n+1}-t^n = T/M$ and $M+1$ nodes so that $t^0=0<t^1<\cdots<t^M=T$.  At the $n$th time node, the state variable $u(x,t)$ is represented by a vector $\bold u ^n = [u_1^n,\cdots,u_N^n]^\top$ for $n = 0,\cdots, M$. For simplicity,~\eqref{2-1} is solved with a conservative first-order upwind scheme with forward Euler for the advection part and center difference with backward Euler for the diffusion part,
\begin{equation}\label{2-2}
u_j^{n+1} = u_j^n-\frac{\Delta t}{\Delta x}(F_{j+1/2}^n-F_{j-1/2}^n)+\frac{\Delta t}{(\Delta x)^2}\left[D_{j+1/2}^{n+1}(u_{j+1}^{n+1}-u_j^{n+1})-D_{j-1/2}^{n+1}(u_j^{n+1}-u_{j-1}^{n+1})\right],
\end{equation}
where 
$$\begin{aligned}
&F_{j+1/2}^n = \frac{F(u_{j+1}^n)+F(u_j^n) }{2} - |a_{j+1/2}^n| \frac{u_{j+1}^n-u_j^n}{2},\\
&a_{j+1/2}^n = \left\{
\begin{aligned}
&\frac{F_{j+1}^n-F_j^n}{u_{j+1}^n-u_j^n}&&\mbox{if}&&u_{j+1}^n\neq u_j^n,\\
&f(u_j)&&\mbox{if}&&u_{j+1}^n=u_j^n,
\end{aligned}
\right.\\
&D_{j+1/2}^{n+1} = \frac{1}{2}(D_j^{n+1}+D_{j+1}^{n+1}).
\end{aligned}$$
In vector form, the above scheme reads 
\begin{equation}\label{2-3}
\bold R(\bold u^{n+1}) = \bold u^{n+1}-\bold u^n+\Delta t(\mathcal D_1^{u} \bold F^n)-\Delta t (\mathcal D_2\bold u^{n+1})=0,
\end{equation}
where %$\odot$ denotes the Hadamard product, 
$\mathcal D_1^u\in\mathbb R^{N\times N}$ and $\mathcal D_2\in\mathbb R^{N\times N}$ are discrete approximations of the first derivative (using upwind) and second derivative (using center difference), respectively. Here $\bold F^n = [F_{1/2}^n, \cdots,F_{N-1/2}^n]^\top$  and $\bold R$ is the vectorized residual of the scheme. Certain CFL condition needs to be satisfied to ensure the stability of the scheme depending on the functional forms of $f$ and $D$. Simulation results obtained with the above method constitute our high-fidelity model (HFM). 

A reduced-order, low-fidelity model (ROM) is constructed from a data set comprising a sequence of solution snapshots collected from the HFM. Let $\bold X$ denote the data matrix, consisting of $m$ snapshots of $\mathbf u$,
\begin{equation}\label{2-4}
\bold X = \begin{bmatrix}
|&|&&|\\
\bold u^1&\bold u^2&\cdots&\bold u^{m}\\
|&|&&|
\end{bmatrix}, \qquad \mathbf X \in \mathbb R^{N \times m}.
\end{equation}
Two alternative strategies for building a ROM from these data, both grounded in Singular Value Decomposition (SVD), are described below.

\subsection{POD}
\label{sec:convPOD}

The conventional POD method generates a ROM by using a low-dimensional basis extracted from the data $\mathbf X$ in~\eqref{2-4} to project the dynamics $\mathbf u(t)$ onto a lower-dimensional hyperplane. If the data matrix $\mathbf X \in \mathbb R^{N \times m}$ has rank $K \leq\min\{N,m\}$, then the POD modes are constructed by using a reduced SVD,
\begin{equation}\label{2-5}
\bold X = \bold U\bold \Sigma\bold V^*,
\end{equation}
where $\bold U \in \mathbb C^{N\times K}$ is the matrix of $K$ orthonormal columns of length $N$; $\bold \Sigma = \mathbb R^{K\times K}$ is the diagonal matrix with real diagonal elements $\sigma_1 \ge \sigma_2 \ge \cdots \ge \sigma_K > 0$; $\bold V \in \mathbb C^{m\times K}$ is the matrix of $K$ orthonormal columns of length $m$; and the superscript $^*$ denotes its conjugate transpose. A reduced-order model is constructed by choosing a rank $r$ ($r \ll K$), which satisfies the energy criteria
\begin{equation}\label{2-6}
r = \min_{k} \left\{\frac{\sigma_k}{\sum_{k=1}^K\sigma_k}<\varepsilon\right\},
\end{equation}
where $\varepsilon$ is a user-supplied small number ($\varepsilon = 10^{-4}$ in all our numerical examples). Next, the matrix $\bold U \in \mathbb C^{N\times K}$ is replaced with a matrix $\boldsymbol \Phi \in \mathbb C^{N \times r}$ comprising $r$ orthogonal columns of length $N$, 
\begin{equation}\label{2-7}
\boldsymbol \Phi = \begin{bmatrix}
|&|&&|\\
\boldsymbol \phi_1&\boldsymbol \phi_2&\cdots& \boldsymbol\phi_r\\
|&|&&|
\end{bmatrix}. 
\end{equation} 
The orthonormal vectors $\{\boldsymbol\phi_1,\cdots,\boldsymbol\phi_r\}$ form a POD basis. Finally, a ROM (low-fidelity solution) is constructed by the Galerkin projection of $\mathbf u$ onto the low-dimensional space spanned by the POD basis,
\begin{equation}\label{2-8}
\bold u_\text{POD}^{n+1} = \sum_{k = 1}^r\hat u_k^{n+1} \boldsymbol\phi_k = \boldsymbol \Phi\hat{\bold u}^{n+1}.
\end{equation}
Substituting (\ref{2-8}) into (\ref{2-3}) and projecting onto the low-dimensional space, yields an equations for the vector of coefficients $\hat{\bold u}^{n+1}$:
\begin{equation}\label{2-9}
\boldsymbol \Phi^\top\bold R( \boldsymbol \Phi\hat{\bold u}^{n+1})=0.
\end{equation}

To deal with the nonlinearity of~\eqref{2-9} numerically, one might use Newton iteration or other efficient methods~\cite{barrault2004empirical, chaturantabut2010nonlinear}.

\subsection{DMD}

We start by recasting the spatially discretized~\eqref{2-1} in the form of a general nonlinear dynamic system 
\begin{equation}\label{3-7}
\frac{\text d\bold u}{\text dt} = \mathcal N(\bold u),
\end{equation}
where $\bold u \in \mathcal M \subset \mathbb R^N$, $\mathcal M$ is a smooth $N$-dimensional manifold, and $\mathcal N$ is a finite dimensional nonlinear operator. Given a flow map $\mathcal N_t :\mathcal M \to \mathcal M$,
\begin{equation}\label{3-8}
\mathcal N_t (\bold u(t_0)) = \bold u(t_0+t) = \bold u(t_0)+\int_{t_0}^{t_0+t} \mathcal N(\bold u(\tau))\bold d \tau,
\end{equation}
the time-discretized counterpart of~\eqref{3-7} is 
\begin{equation}\label{3-9}
\bold u^{n+1} = \mathcal N_t(\bold u^n).
\end{equation}

The DMD method approximates the modes of the so-called \textit{Koopman operartor}:
\begin{definition}[Koopman operator \cite{kutz2016dynamic}]\label{def_koopman}
The Koopman operator $\mathcal K$ for the nonlinear dynamic system~\eqref{3-7} is an infinite-dimensional linear operator that acts on all observable functions $g: \mathcal M\to \mathbb R$ such that
\begin{equation}\label{3-10}
\mathcal K g(\bold u) = g(\mathcal N(\bold u)).
\end{equation}
The discrete-time Koopman operator $\mathcal K_t$ for the discrete dynamic system (\ref{3-9}) is defined by
\begin{equation}\label{3-11}
\mathcal K_t g(\bold u^{n}) = g(\mathcal N_t(\bold u^n)) =g(\bold u^{n+1}).
\end{equation}
\end{definition}

In practice, the most accessible observable is usually the state itself. Thus conventional DMD method  generates a ROM by seeking a truncated finite approximation of the Koopman operator $\mathcal K$ coorperating with the chosen observable function $g$ as an identity map. This is done by splitting the data matrix $\bold X$ into two,
\begin{equation}\label{2-10}
\bold X_1 = \begin{bmatrix}
|&|&&|\\
\bold u^1&\bold u^2&\cdots&\bold u^{m-1}\\
|&|&&|
\end{bmatrix}, \qquad 
\bold X_2 = \begin{bmatrix}
|&|&&|\\
\bold u^2&\bold u^3&\cdots&\bold u^{m}\\
|&|&&|
\end{bmatrix},
\end{equation}
and using these two datasets to approximate the eigenvalues and eigenvectors of $\mathcal K$ by means of the following algorithm~\cite{kutz2016dynamic}.
\begin{algo}{DMD algorithm}
\begin{enumerate}
\item Compute SVD of matrix $\bold X_1 \approx \bold U\boldsymbol \Sigma \bold V^*$ with $\bold U \in \mathbb C^{N\times r}$, $\boldsymbol \Sigma \in \mathbb R^{r\times r}$ and $\bold V\in \mathbb C^{r\times m}$, where $r$ is the truncated rank chosen by a certain criteria, e.g. (\ref{2-6}).
\item Compute $\tilde{\bold K}=\bold U^*\bold X_2\bold V\boldsymbol\Sigma^{-1}$ as an $r\times r$ low-rank approximation of $\mathcal K$.
\item Compute eigendecomposition of $\tilde{\bold K}$: $\tilde{\bold K} \bold W = \bold W\boldsymbol\Lambda$, $\boldsymbol \Lambda = (\lambda_k)$.
\item Reconstruct eigendecomposition of $\mathcal K$, whose eigenvalues and eigenvectors are $\boldsymbol\Lambda$ and $\boldsymbol \Phi  = \bold U\bold W$, respectively.
\end{enumerate}
\end{algo}
Each column of $\boldsymbol \Phi$ is a DMD mode corresponding to a particular eigenvalue in $\boldsymbol\Lambda$. With the approximated eigenvalues and eigenvectors of $\mathcal K$ in hand, the projected future solution can be constructed analytically for all times in the future. In particular, at each future time $t = t^{n+1}$,
\begin{equation}\label{2-11}
\bold u_\text{DMD}^{n+1} =\boldsymbol\Phi\boldsymbol\Lambda^{n+1}\bold b,  \ \ n>m,
\end{equation}
where $\bold b = \Phi^{-1}\bold u^{0}$ is an $r\times1$ vector representing the initial amplitude of each mode. Note that no more iterations are needed in the prediction. The solution at any future time can be approximated directly with~\eqref{2-11} using only information from the first $m$ snapshots of the HFM.

\subsection{Challenge Posed by Translational Problems}

Both POD and DMD have been used to construct LFMs for a wide range of problems with high accuracy. However, ROMs constructed with such SVD-based methods are known to have poor performance for translational problems, such as an advection-dominated version of~\eqref{2-1}. To illustrate this phenomenon, we consider a linear advection-diffusion equation, i.e.,~\eqref{2-1} with constant $f$ and $D$, defined on $(x,t) = [0,2] \times [0,1]$. This equation is subject to the initial condition $u(x,0) = 0.5 \exp[-(x-0.3)^2/0.05^2]$ and the boundary conditions $u(0,t) = u(2,t) = 0$. The space domain $[0,2]$ is discretized into $N=2000$ intervals and time domain $[0,1]$ is discretized into $M =1000$ steps. Both DMD and POD algorithms use the same dataset consisting of $m = 250$ snapshots. 

\begin{figure}[htbp]
\centering
\includegraphics{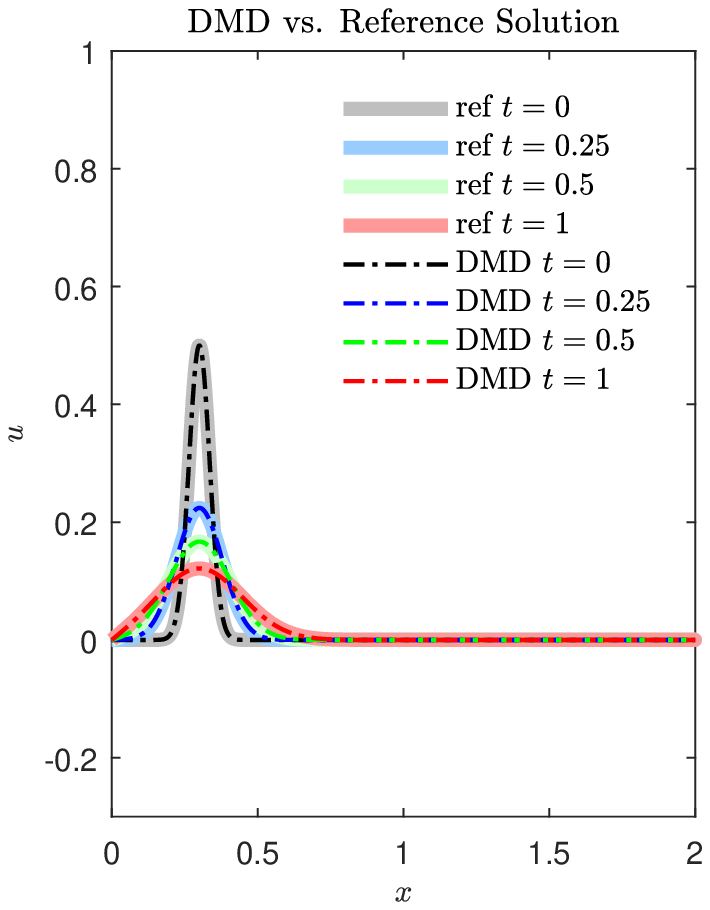}
\includegraphics{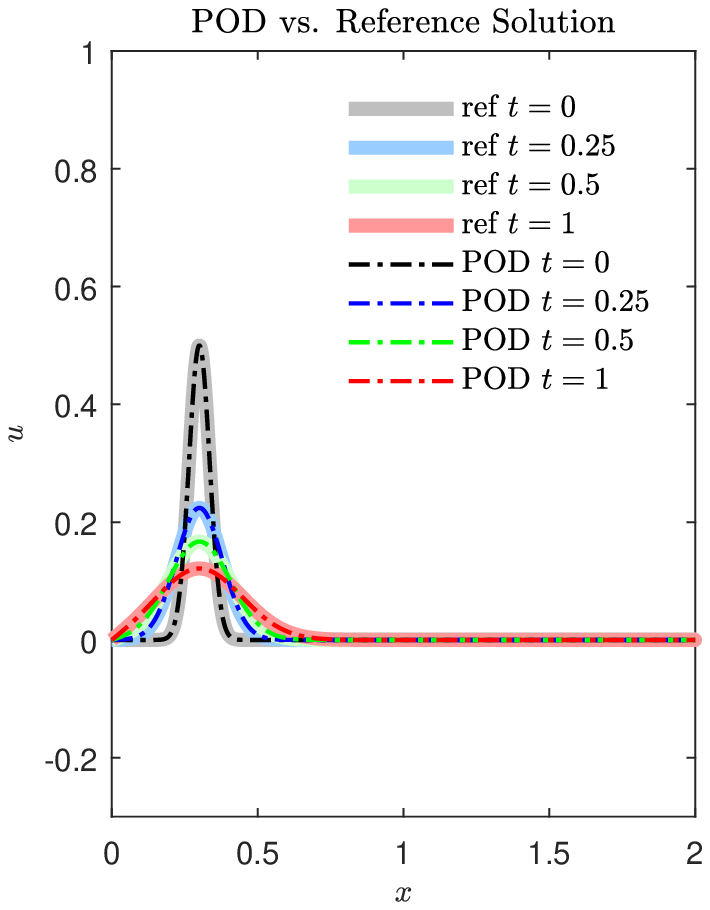}
\caption{Solution profiles $u(x,\cdot)$, for several times $t$, in the diffusion-dominated regime. These profiles are computed with DMD (left) and POD (right), and compared with the reference solution. \label{fig:diff-dom} }
\end{figure}

To achieve a diffusion-dominated regime, we set $f = 10^{-4}$ and $D = 10^{-2}$ in some consistent units. Figure~\ref{fig:diff-dom} provides a visual comparison of the reference solution with its counterparts obtained with DMD and POD, both with $r = 20$ SVD rank truncation. Although not shown here, and consistent with the earlier findings reported in~\cite{lu2019error}, the DMD- and POD-based LFMs are of spectral accuracy in a relatively small subspace of time ($t < 0.3$), with the relative error increasing with time. POD has slightly better accuracy than DMD due to the iterations in the subspace, but DMD is considerably faster because of its iteration-free nature.

\begin{figure}[htbp]
\centering
\includegraphics{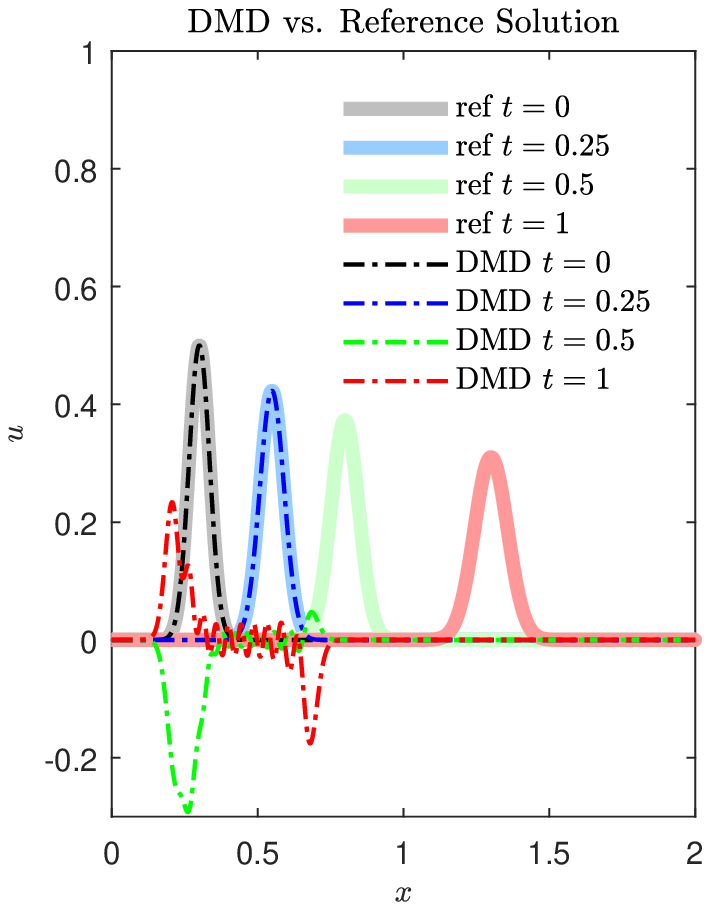}
\includegraphics{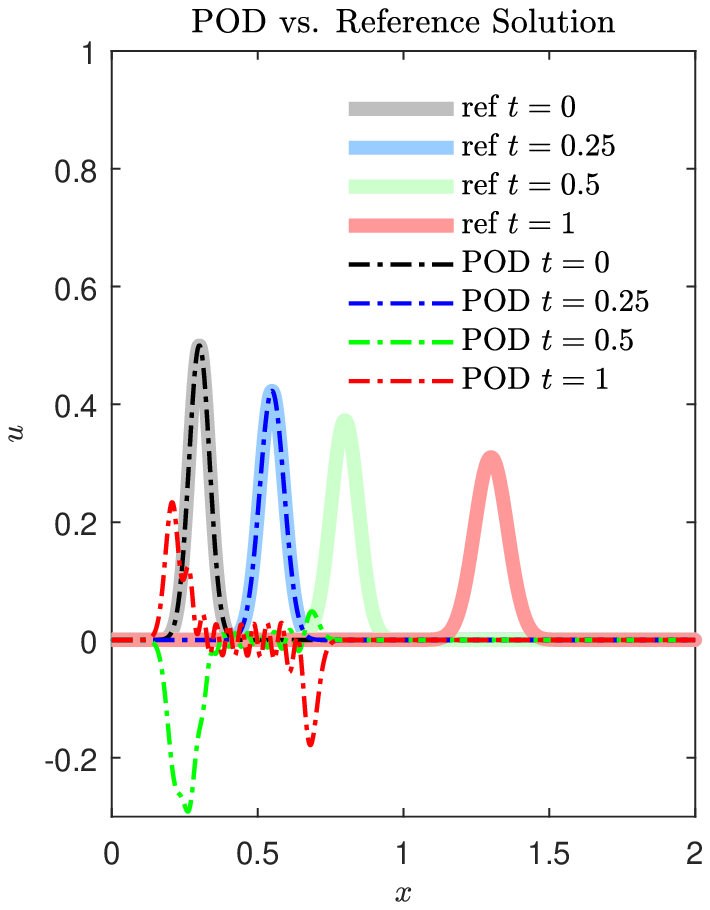}
\includegraphics{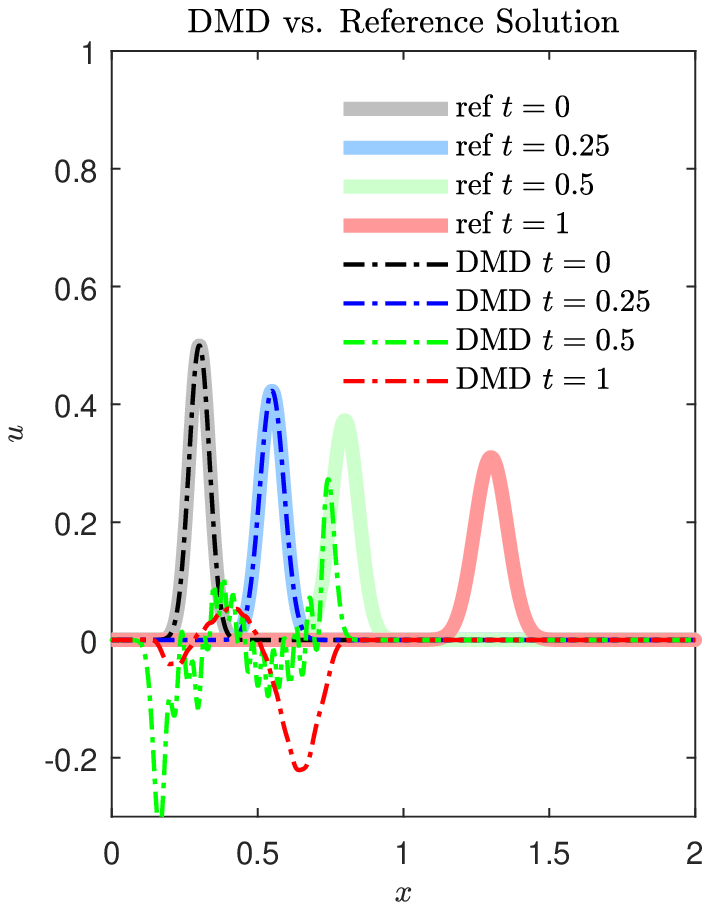}
\includegraphics{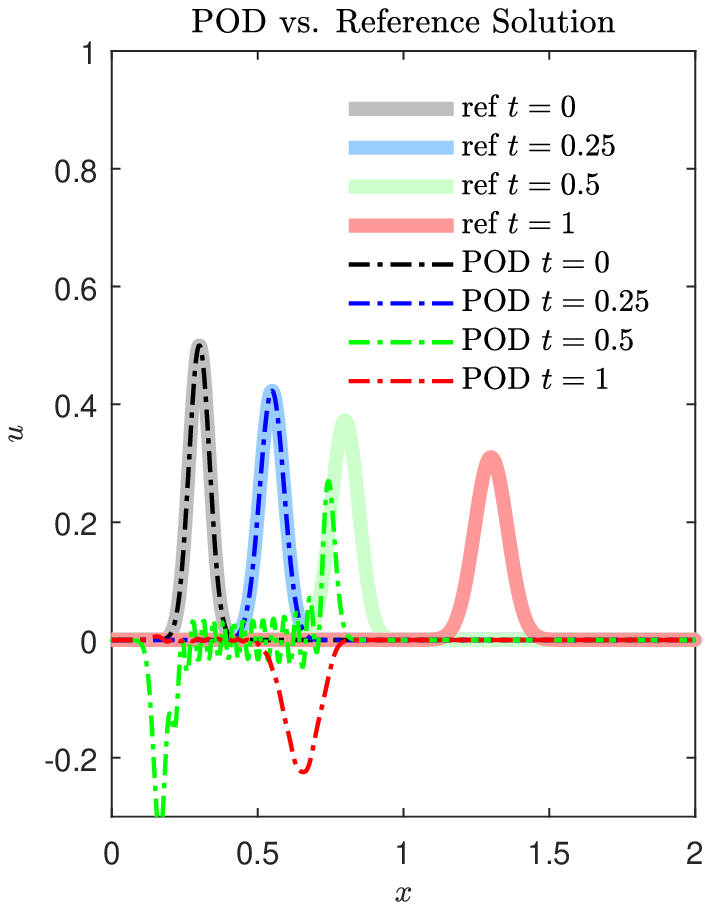}
\caption{Solution profiles $u(x,\cdot)$, for several times $t$, in the advection-dominated regime. These profiles are computed with DMD (left column) and POD (right column) using the SVD rank of $r = 20$ (top row) and $r=30$ (bottom row), and compared with the reference solution. \label{fig:adv-dom} }
\end{figure}

\begin{figure}[htbp]
\centering
\includegraphics{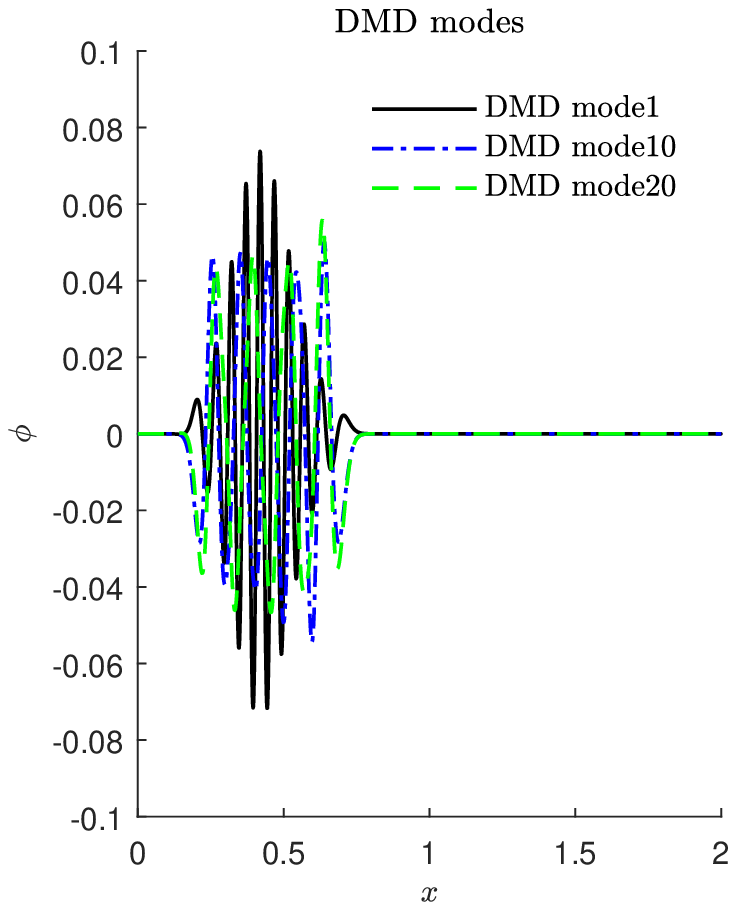}
\includegraphics{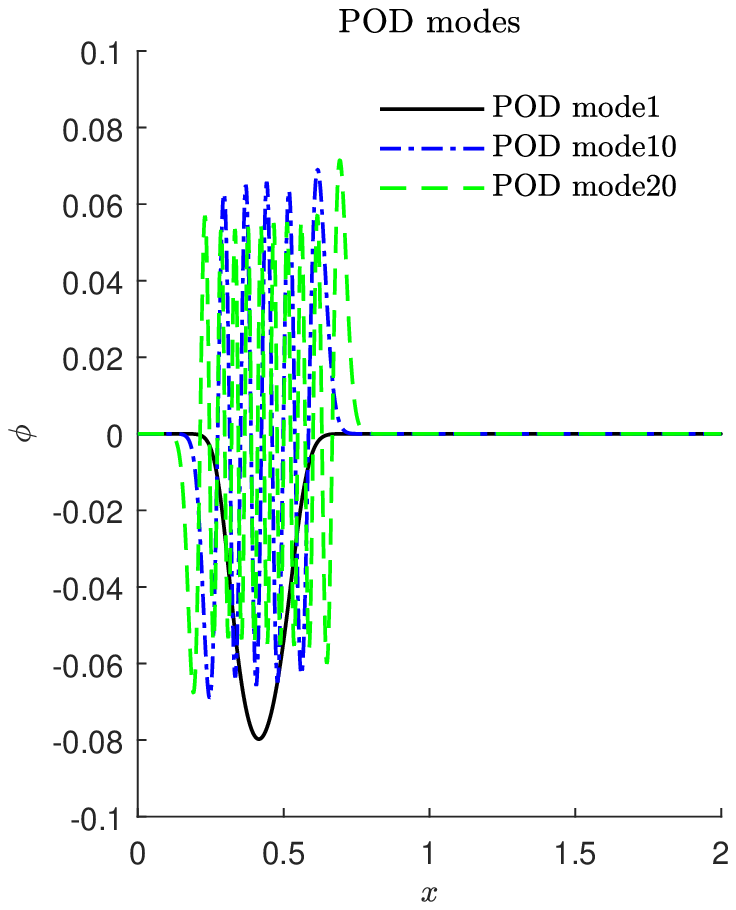}
\caption{Three of the dominant DMD modes (left column) and POD modes (right column) extracted from the first 250 snapshots. \label{fig:modes}}
\end{figure}

An advection-dominated regime is achieved by setting $f = 1.0$ and $D = 10^{-3}$. Figure~\ref{fig:adv-dom} reveals that both DMD and POD fail to capture the system dynamics, yielding unphysical (oscillatory and negative) predictions. This failure cannot be remedied by increasing the SVD rank truncation $r$: increasing $r$ from 20 to 30 does not improve the prediction's accuracy, either quantitatively or qualitatively. These results highlight the main challenge translational problems pose for the SVD-based methods.  Given the first 250 snapshots of the high-fidelity solution, SVD extracts dominant DMD/ POD modes $\phi_i$ from the region the wave has encountered; in our example, the subdomain $[0,1]$. As time increases, the wave solution encounters other parts of the computational domain; in our example, at later times the dominant signal lies mostly in the subdomain $[1,2]$. Specifictly,  one can observe that the dominant DMD/POD modes have fluctuations only  in the subdomain $[0,1]$ and stay flat $0$ in the subdomain $[1,2]$ in Figure~\ref{fig:modes}. It is therefore not surprising that a ROM constructed from dominant modes in $[0,1]$ does not serve as an accurate surrogate for the rest of the computational domain.

%%%%%%%%%%%%%%%%%%%%%%%%%%%%%%
\section{Lagrangian Reduced-Order Models}
\label{sec:Lagrange}
%%%%%%%%%%%%%%%%%%%%%%%%%%%%%%

Motivated by construction of a POD-based ROM for the advection-diffusion equation~\eqref{2-1} within the Lagrangian framework~\cite{mojgani2017lagrangian}, we propose a Lagrangian DMD. In the semi-Lagrangian frame,~\eqref{2-1} is written as
\begin{equation}\label{3-1}
\left\{
\begin{aligned}
&\frac{\text dX(t)}{\text dt} = f(u(X(t),t)),\\
&\frac{\text du(x,t)}{\text dt}\bigg\rvert_{x=X(t)} = \left[\frac{\partial}{\partial x}\left(D(x,t,u)\frac{\partial u(x,t)}{\partial x}\right)\right]\bigg\rvert_{x=X(t)}. 
\end{aligned}
\right.
\end{equation}
Applying a first-order discretization to this system gives
\begin{equation}\label{3-2}
\left\{
\begin{aligned}
&\tilde u_j^n = \mathcal P_0(u_j^n),\\
&\tilde u_j^{n+1} = \tilde u_j^n+\Delta t \frac{D_{j+1/2}^{n+1}(\tilde u_{j+1}^{n+1}-\tilde u_j^{n+1})-D_{j-1/2}^{n+1}(\tilde u_j^{n+1}-\tilde u_{j-1}^{n+1})}{(\Delta x)^2},\\
&u_j^{n+1} = \mathcal P_{n} (\tilde u_j^{n+1}),\\
&x_j^{n+1} = x_j^n+\frac{\Delta t}{2}(f(u_j^n)+ f(u_j^{n+1})),
\end{aligned}
\right.
\end{equation}
where $\mathcal P_n$ stands for the interpolation in the grid $\bold x^n = (x_1^n,\cdots, x_N^n)^\top$ and $\bold x^0$ is the starting uniform grid.

%\begin{equation}\label{3-2}
%\left\{
%\begin{aligned}
%&\tilde u_j^n = \mathcal P_0(u_j^n),\\
%&\tilde u_j^{n+1} = \tilde u_j^n+\Delta t \frac{D_{j+1/2}^{n+1}(\tilde u_{j+1}^{n+1}-\tilde u_j^{n+1})-D_{j-1/2}^{n+1}(\tilde u_j^{n+1}-\tilde u_{j-1}^{n+1})}{(\Delta x)^2},\\
%& u_j^{n+1} = \mathcal P_n(\tilde u_j^{n+1}),\\
%&
%\end{aligned}
%\right.
%\end{equation}

Or, in vector form,
%In practice, one would compute $\bold u^{n+1}$ first and approximate $\bold x^{n+1}$ using mid-point rule:
%\begin{equation}
%\left\{
%\begin{aligned}
%&\bold u^{n+1} = \bold u^n+\Delta t \bold g (\bold x^{n+1},\bold u^{n+1})\odot \mathcal D_2\bold u^{n+1},\\
%&\bold x^* = \bold x^n+\frac{\Delta t}{2}\bold f(\bold x^n,\bold u^n),\\
%&\bold x^{n+1} = \bold x^n+\frac{\Delta t}{2} \left[\bold f(\bold x^{*},\bold u^{n}(\bold x^*))+\bold f(\bold x^*,\bold u^{n+1}(\bold x^*))\right].
%\end{aligned}
%\right.
%\end{equation}
\begin{equation}\label{3-3}
\left\{
\begin{aligned}
&\bold R_x(\bold x^{n+1}) \equiv \bold x^{n+1}-\bold x^n-\frac{\Delta t}{2}(\bold f(\bold u ^{n})+\bold f(\bold u ^{n+1}))=0,\\
&\bold R_u(\bold {\tilde u}^{n+1}) \equiv \bold {\tilde u}^{n+1}-\bold {\tilde u}^n-\Delta t\mathcal D_2\bold {\tilde  u}^{n+1}=0.
\end{aligned}
\right.
\end{equation}
Here $\bold x^n = [x_1^n,\cdots, x_N^n]^\top$ denotes the locations of the Lagrangian computational grid at the $n$th time step, $\bold{\tilde u}^n$ is the interpolation from the Lagrangian grid to the Eulerian grid, and $\mathcal D_2$ represents the discrete approximation of the second derivative on the uniform Eulerian grid at the $n$th time step.

\subsection{POD}

We arrange $m$ snapshots of the HFM in the Lagrangian framework into a data matrix of size $2N\times m$,
\begin{equation}\label{3-4}
\bold X = \begin{bmatrix}
|&|&&|\\
\bold x^1&\bold x^2&\cdots&\bold x^{m}\\
|&|&&|\\
\bold u^1&\bold u^2&\cdots&\bold u^{m}\\
|&|&&|
\end{bmatrix}.
\end{equation}
Applying the conventional POD of section~\ref{sec:convPOD} to the data matrix in~\eqref{3-4}, one obtains a POD basis $\boldsymbol \Phi$ analogous to~\ref{2-7} for the space-solution vector $[\mathbf x; \mathbf u]^\top$. Then the Lagrangian solution is approximated by 
\begin{equation}\label{3-5}
\begin{bmatrix}
|\\
\bold x_\text{POD}^{n+1}\\
|\\
\bold u_\text{POD}^{n+1}\\
|
\end{bmatrix} = \boldsymbol \Phi 
\begin{bmatrix}
|\\
\bold {\hat x}^{n+1}\\
|\\
\bold {\hat u}^{n+1}\\
|
\end{bmatrix}.
\end{equation}
Inserting (\ref{3-5}) into (\ref{3-3}) and projecting onto the subspace spanned by $\boldsymbol \Phi$, one obtains the solution vector $[\hat{\bold x}^{n+1};\hat{\bold u}^{n+1}]^\top$ by solving the following equation:
\begin{equation}\label{3-6}
\boldsymbol \Phi^T \bold R\left(
\boldsymbol \Phi
\begin{bmatrix}
|\\
\bold {\hat x}^{n+1}\\
|\\
\bold {\hat u}^{n+1}\\
|
\end{bmatrix}
\right)=0.
\end{equation}

Several complications can arise when applying Lagrangian POD in practice. If only Eulerian HFM data are available, i.e., in the absence of the grid deformation $\bold x^n$ computed with an Eulerian HFM, one can construct an optimal Lagrangian basis by following the strategy proposed in~\cite[Sec.~3.3]{mojgani2017lagrangian}. Another potential complication is a Lagrangian grid entanglement. There is no guarantee that an approximation of the Lagrangian moving grid in the low-dimensional subpspace preserves the topological properties of the original HFM simulation. In many cases, e.g., when characteristic lines intersect each other, the Lagrangian grids in the projected space are severely distorted,  inducing numerical instabilities. One strategy for ameliorating this problem is to solve the diffusion step back to stationary Eulerian grid by interpolation between the Eulerian and Lagrangian grids~\cite[Sec.~3.4]{mojgani2017lagrangian}. This procedure can reduce the method's efficiency and accuracy. 

\subsection{DMD}
The fundamental concept behind the Koopman operator theory is to  transform the finite-dimensional nonlinear problem~\eqref{3-9} in the state space into the infinite-dimensional linear problem \eqref{3-11} in the observable space. Compared to Eulerian framework, Lagrangian framework contains more informative physical quantanties as candidates of the obervables, making Lagrangian DMD to fit intuitively and naturally into the Koopman operator theory. We briefly review the ethics in choosing observable functions below for the sake of completeness, followed by a description of our approach for approximating the underlying Koopman operator with Lagrangian DMD. 

Since $\mathcal K_t$ in Definition~\ref{def_koopman} is an infinite-dimensional linear operator, it is equipped with infinitely many eigenvalues $\{\lambda_k\}_{k=1}^{\infty}$ and eigenfunctions $\{\phi_k\}_{k=1}^\infty$. In practice, one deals with a finite number of the eigenvalues and eigenfunctions. The following assumption underpins the finite approximation and is essential to the choice of observables.

\begin{assump}
Let $\mathbf y$ denote a vector of observables,
\begin{equation}\label{3-12}
\bold y^n = \bold g(\bold u^{n}) = \begin{bmatrix}
g_1(\bold u^n)\\
g_2(\bold u^n)\\
\vdots\\
g_p(\bold u^n)
\end{bmatrix}, \qquad g_j: \mathcal M\to \mathbb C \ \mbox{are observable function ($j =1,2,\cdots, p$)}. 
\end{equation} 
If the chosen observable $\bold g$ is restricted to an invariant subspace spanned by eigenfunctions of the Koopman operator $\mathcal K_t$, then it induces a linear operator $\bold K$ that is finite-dimensional and advances these eigen-observable functions on this subspace~\cite{brunton2016koopman}.
\end{assump}

Based on the above assumption, DMD algorithm~\cite{kutz2016dynamic} is applied to approximate the eigenvalues and eigenfunctions of $\bold K$ using the collected snapshot data in the observable space: 

\begin{algo}{Physics-aware DMD algorithm:}
\label{DMDalg}
\begin{itemize}
\item[0.] Create data matrices of observables $\mathbf Y_1$ and $\mathbf Y_2$ as
\begin{equation}\label{3-13}
\bold Y_1 = \begin{bmatrix}
|&|&&|\\
\bold y^1&\bold y^2&\cdots&\bold y^{m-1}\\
|&|&&|
\end{bmatrix},\bold Y_2 = \begin{bmatrix}
|&|&&|\\
\bold y^2&\bold y^3&\cdots&\bold y^{m}\\
|&|&&|
\end{bmatrix},
\end{equation}
where each column is given by $\bold y^k = \bold g(\bold u^k)$.
\item[1.] Compute SVD of the matrix $\bold Y_1 \approx \bold U\boldsymbol \Sigma \bold V^*$ with $\bold U \in \mathbb C^{p\times r}, \boldsymbol\Sigma \in \mathbb R^{r\times r}, \bold V\in \mathbb C^{r\times m}$, where $r$ is the truncated rank chosen by certain criteria.
\item [2.] Compute $\tilde{\bold K}=\bold U^*\bold X'\bold V\boldsymbol\Sigma^{-1}$ as an $r\times r$ low-rank approximation for $\bold K$.
\item [3.] Compute eigendecomposition of $\tilde{\bold K}$: $\tilde{\bold K} \bold W = \bold W\boldsymbol\Lambda$, $\boldsymbol\Lambda = (\lambda_k)$.
\item [4.] Reconstruct eigendecomposition of $\bold K$. Eigenvalues are $\boldsymbol\Lambda$ and eigenvectors are $\boldsymbol\Phi  = \bold U\bold W$.
\item [5.] Future $\bold y_{\text{DMD}}^{n+1}$ can be predicted by
\begin{equation}\label{3-14}
\bold y_{\text{DMD}}^{n+1} =\boldsymbol \Phi\boldsymbol\Lambda^{n+1}\bold b, \ \ n>m
\end{equation}
with $\bold b =\boldsymbol \Phi^{-1}\bold y_1$.
\item [6.] Transform from observables back to state-space:
\begin{equation}\label{3-15}
\bold u_\text{DMD}^n =\bold g^{-1}(\bold y_\text{DMD}^n).
\end{equation}
\end{itemize}
\end{algo}

In data-driven modeling, judicious selection of the observables is crucial to the accuracy and efficiency of a Koopman operator's approximation. Identification of  general rules for choosing the observables continues to be a subject of ongoing research. For example, the use of measurements of the state variable $u(x,t)$ as an observable led to the poor DMD performance in the advection-dominated regime (Figure~\ref{fig:adv-dom}). %It is worth to mention that several recent work has been done to approximate the Koopman operator using machine learning and deep learning \cite{takeishi2017learning,yeung2017learning,morton2018deep}. 
A Lagrangian formulation of the problem provides a means of identification of optimal observables. Indeed, the physics of advection-dominated systems suggests that the location of a moving particle is a key quantity, which is as important as the value of the state variable at that location. It is therefore natural to introduce an observable function that keeps track of both. Thus we choose our observable to be
\begin{equation}
\bold y^n = \bold g(\bold u^n) = \begin{bmatrix}
g_1^n\\
g_2^n
\end{bmatrix}, \quad \mbox{with} \ g_1 = \bold x^n, \ g_2 = \bold u^n.
\end{equation}
Then, we follow Algorithm~\ref{DMDalg}.

%%%%%%%%%%%%%%%%%%%%%%%%%%
\section{Numerical Experiments}
\label{sec:example}
%%%%%%%%%%%%%%%%%%%%%%%%%%

To ascertain the accuracy and robustness of the Lagrangian DMD, we use it to construct ROMs for a series of linear and nonlinear advection-dominated problems. In all tests, the reference solutions are computed in the Eulerian framework using~\eqref{2-2}. The space domain $[0,2]$ is discretized into $N = 2000$ intervals and the time domain $[0,1]$ is discretized into $M = 1000$ steps. Both Lagrangian DMD and Lagrangian POD algorithms use $m = 250$  snapshots (up to $t=0.25$) as a training dataset. The rank truncation criteria~\eqref{2-6} with $\varepsilon = 10^{-8}$ is used. The error bound derived in \cite{lu2019error} is reviewed in the following theorem and plotted in each example as an estimate of the observable.

\begin{theo}
Define the global truncation error
\begin{equation}\label{error}
\bold e^n = \bold y^n -\bold y_{DMD}^n.
\end{equation}
Then, for $n\geq m$,
\begin{equation}
\mathcal E^n =\|\bold e^n\|_2\leq \|\boldsymbol \Phi^{-1}\|_F[\|\bold e^m\|_2+(n-m)\varepsilon_m],
\end{equation}
where $\varepsilon_m$ is a constant only depending on the number of snapshots $m$.
\end{theo}

\subsection{Linear Advection Equation}

We start by considering~\eqref{2-1} with $f \equiv 1$ and $D \equiv 0$. The resulting linear advection equation is defined on $(x,t) \in (0,2) \times (0,1]$, and is subject to the initial condition
\begin{subequations}\label{4-1}
\begin{equation}
%\left\{
%\begin{aligned}
%&\frac{\partial u(x,t)}{\partial t}+\frac{\partial u(x,t)}{\partial x} = 0, ,\\
u(x,t=0) = u_0(x) \equiv \frac{1}{2} \exp\left[- \left(\frac{x-0.3}{0.05} \right)^{\!\!2} \right]
\end{equation}
and boundary conditions
\begin{equation}
u(0,t) = u(2,t) = 0.
%\end{aligned}
%\right.
\end{equation}
\end{subequations}

\begin{figure}[htbp]
\centering
\includegraphics{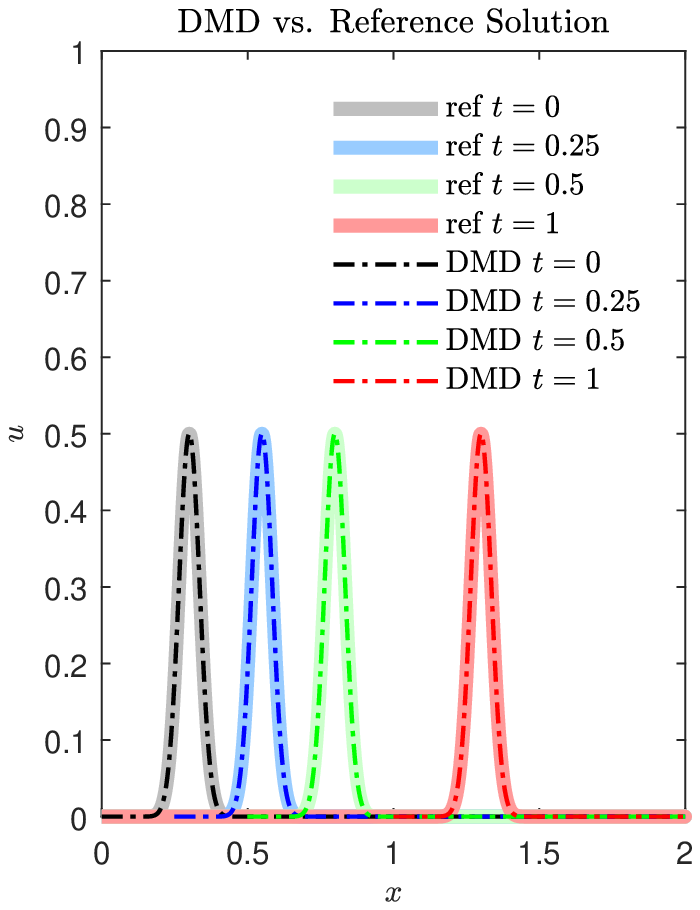}
\includegraphics{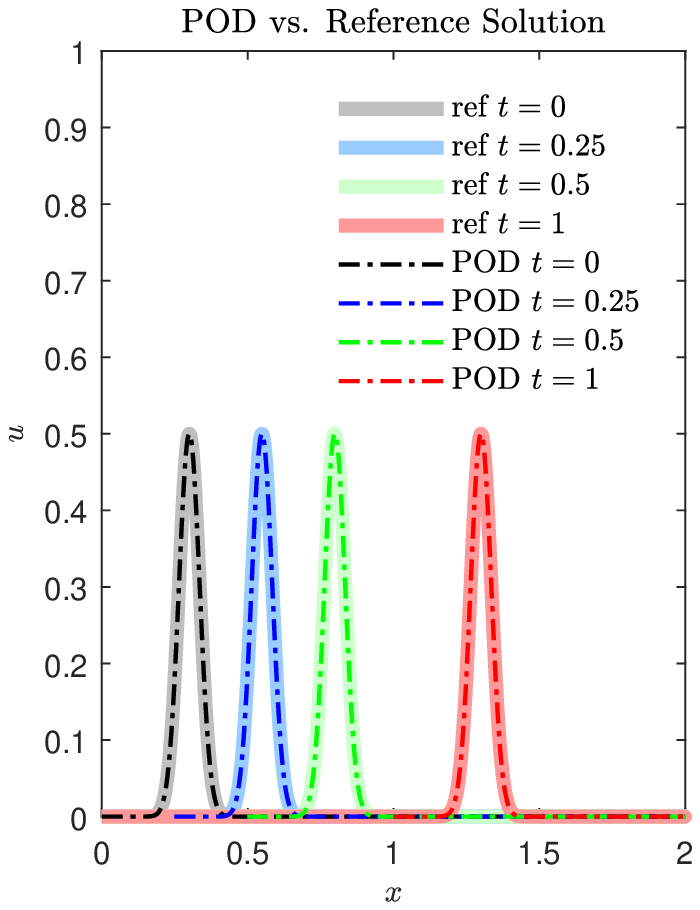}
\caption{Solutions of the linear advection equation, $u(x,t)$, alternatively obtained with the numerical method~\eqref{2-2} and the ROMs constructed via Lagrangian DMD (left) and Lagrangian POD (right). \label{fig:test1-vis} }
\end{figure}

Figure~\ref{fig:test1-vis} provides a visual comparison between the reference solution, obtained with the numerical scheme~\eqref{2-2}, and solutions of the ROMs constructed with either Lagrangian DMD or Lagrangian POD. Unlike their conventional (Eulerian) counterparts (see Figure~\ref{fig:adv-dom}), both Lagrangian DMD and Lagrangian POD capture the solution dynamics in the extrapolating mode, i.e., for $t > 0.25$.

\begin{figure}[htbp]
\centering
\includegraphics{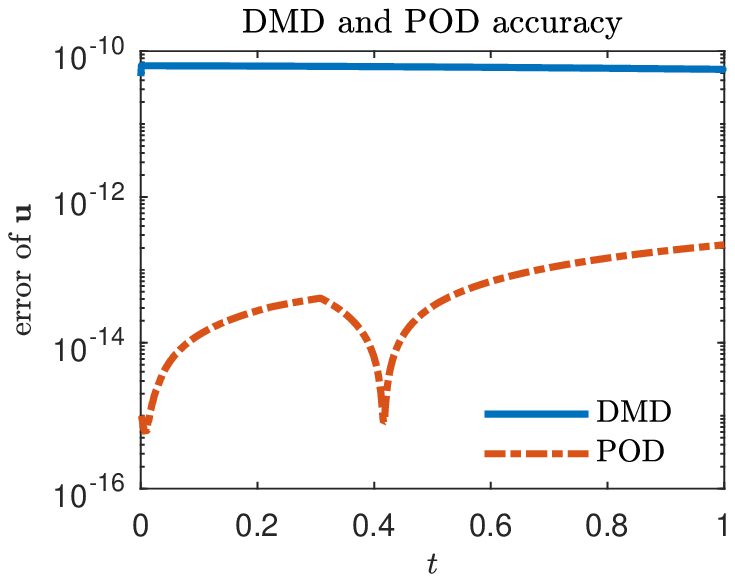}
\includegraphics{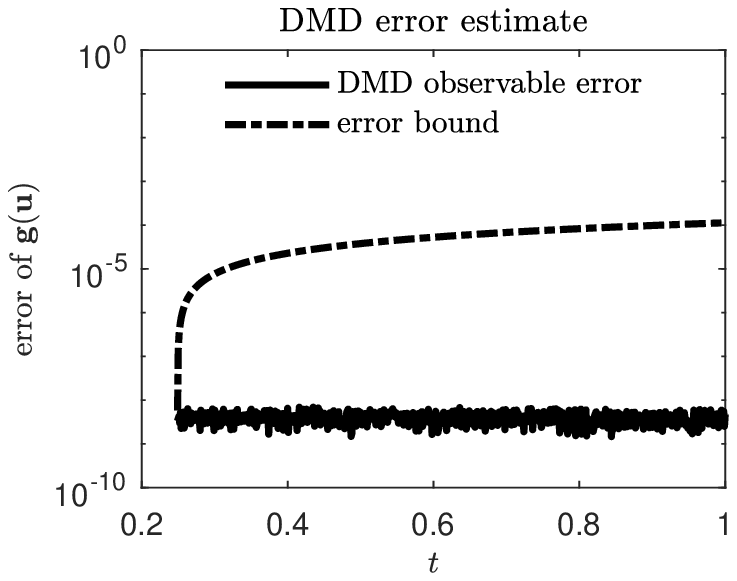}
\caption{Errors $\mathcal E$ of the Lagrangian DMD- and POD-based ROMs for the linear advection equation: error in reconstructing the state variable $u(x,t)$ (left) and its observables $\mathbf g(u)$ (right). The error bound for $\mathbf g(u)$ is derived in~\cite{lu2019error}. \label{fig:test1-err} }
\end{figure}

A more quantitative comparison of the relative performance of the two SVD-based strategies is presented in Figure~\ref{fig:test1-err} in terms of the global truncation error defined in~\eqref{error}. Both Lagrangian DMD and Lagrangian POD capture the advection with high accuracy. Due to the linearity and conservation property of this problem, the ROMs constructed by the two methods are of machine error. Thus, the error bound developed in~\cite{lu2019error} is not tight but sufficient to serve as an indicator of successful approximation.

\subsection{Linear Advection-Diffusion Equation}

Next, we consider~\eqref{2-1} with $f \equiv 1$ and $D \equiv 0.01$. The resulting linear advection-diffusion equation is defined on $(x,t) \in (0,2) \times (0,1]$, and is subject to the initial and boundary conditions~\eqref{4-1}. The choice of the parameter values $f$ and $D$ ensures that the system is in the advection-dominated regime, for which the conventional POD and DMD fail.

\begin{figure}[htbp]
\centering
\includegraphics{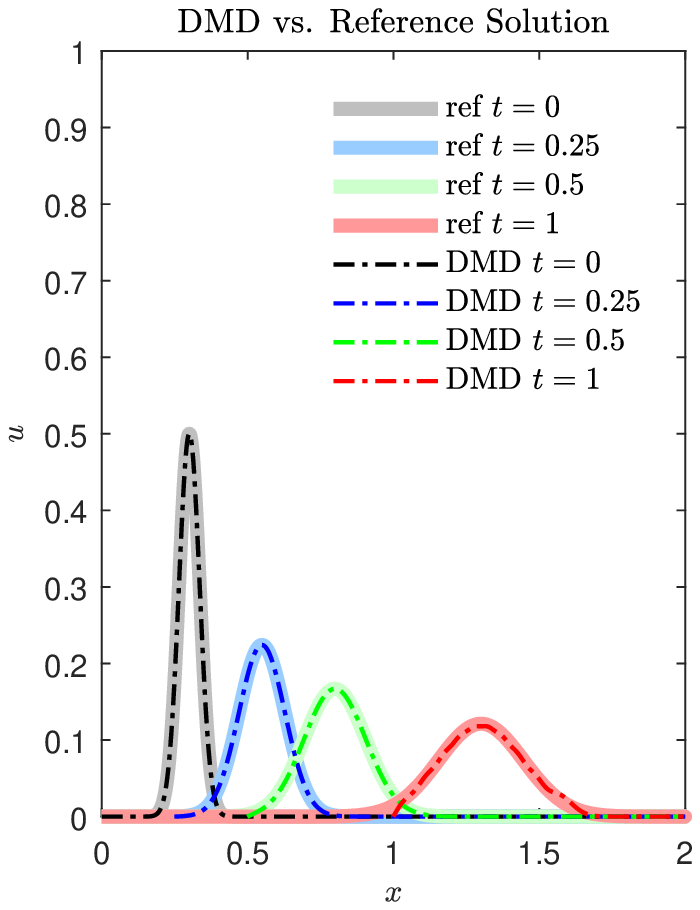}
\includegraphics{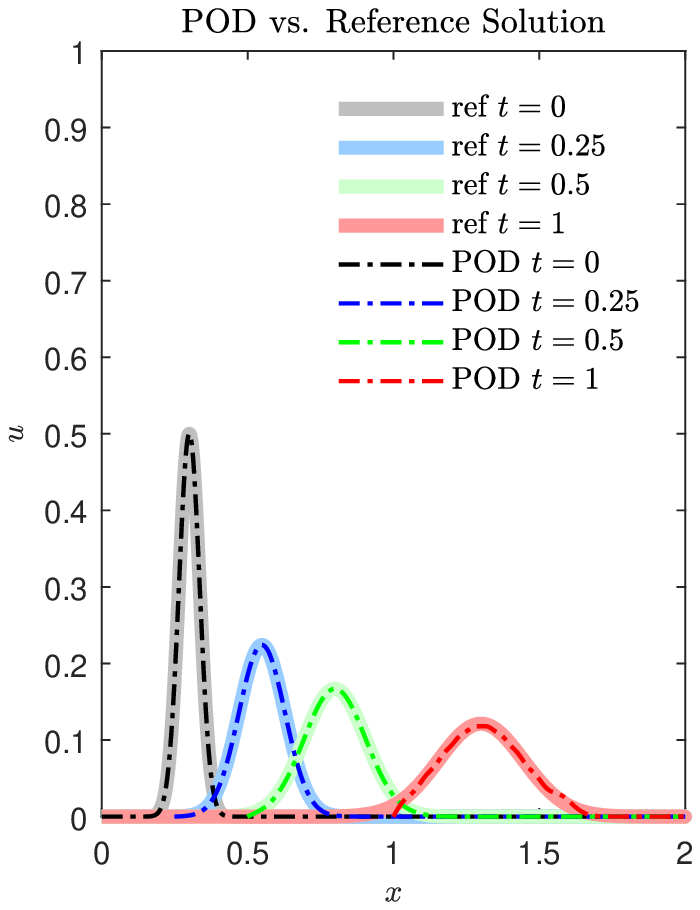}
\caption{Solutions of the linear advection-diffusion equation, $u(x,t)$, alternatively obtained with the numerical method~\eqref{2-2} and the ROMs constructed via Lagrangian DMD (left) and Lagrangian POD (right). \label{fig:test2-vis} }
\end{figure}

Figure~\ref{fig:test2-vis} provides a visual comparison between the reference solution $u(x,t)$ and those obtained with the ROMs. The latter capture the system's dynamics, although their estimates of the solution tails become less accurate with time. This suggests that Lagrangian DMD and POD are capable of identifying the low-rank structure of the advection-diffusion dynamics in the advection-dominated regime.

\begin{figure}[htbp]
\centering
\includegraphics{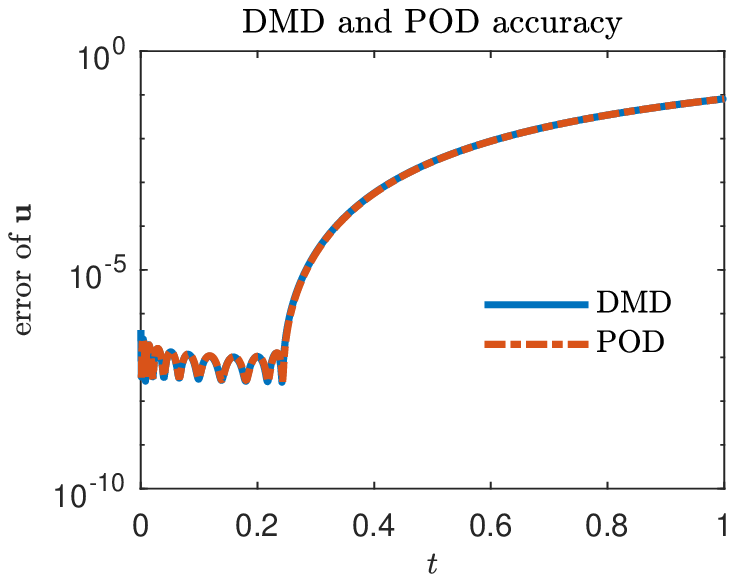}
\includegraphics{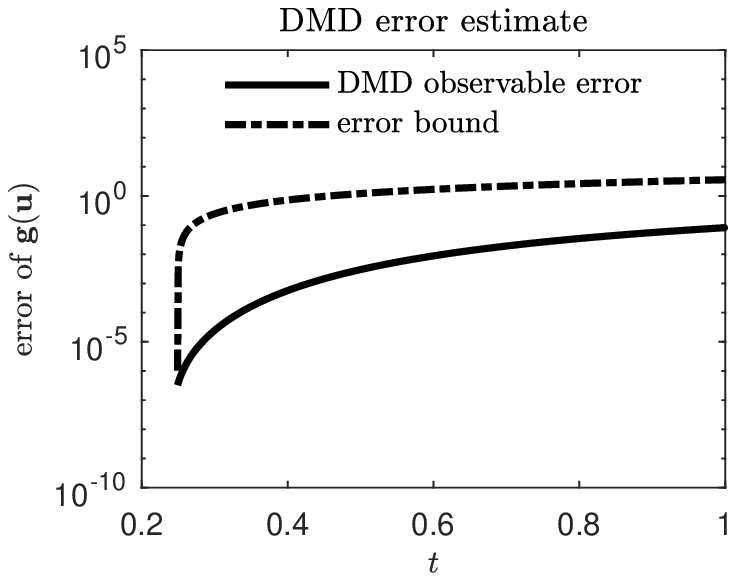}
\caption{Errors $\mathcal E$ of the Lagrangian DMD- and POD-based ROMs for the linear advection equation: error in reconstructing the state variable $u(x,t)$ (left) and its observables $\mathbf g(u)$ (right). The error bound for $\mathbf g(u)$ is derived in~\cite{lu2019error}. \label{fig:test2-err} }
\end{figure}

Figure~\ref{fig:test2-err} indicates that the Lagrangian DMD and POD have a near identical accuracy, which deteriorates with extrapolation time $t > 0.25$. The error bound for the DMD estimate of the observable $g(u)$ is appreciably tighter than in the case of advection (Figure~\ref{fig:test1-err}). With the error bounds, one can design an algorithm combining short-term computation of HFM with long-term computation of LFM.

\subsection{Inviscid Burgers Equation}

The inviscid Burgers equation is recovered from~\eqref{2-1} by setting $f \equiv u$ and $D \equiv 0$. We define this equation on $(x,t) \in (0,2\pi) \times (0,1]$,  subject to the initial conditions
\begin{subequations}\label{4-3}
\begin{equation}
%\left\{
%\begin{aligned}
%&\frac{\partial u(x,t)}{\partial t}+u(x,t)\frac{\partial u(x,t)}{\partial x} = 0, \ [x,t]\in[0,2\pi]\times[0,1],\\
u(x,t=0) = u_0(x) \equiv 1+\sin(x)
\end{equation}
and the periodic boundary conditions
\begin{equation}
u(0,t) = u(2\pi,t).
%\end{aligned}
%\right.
\end{equation}
\end{subequations}

\begin{figure}[htbp]
\centering
\includegraphics{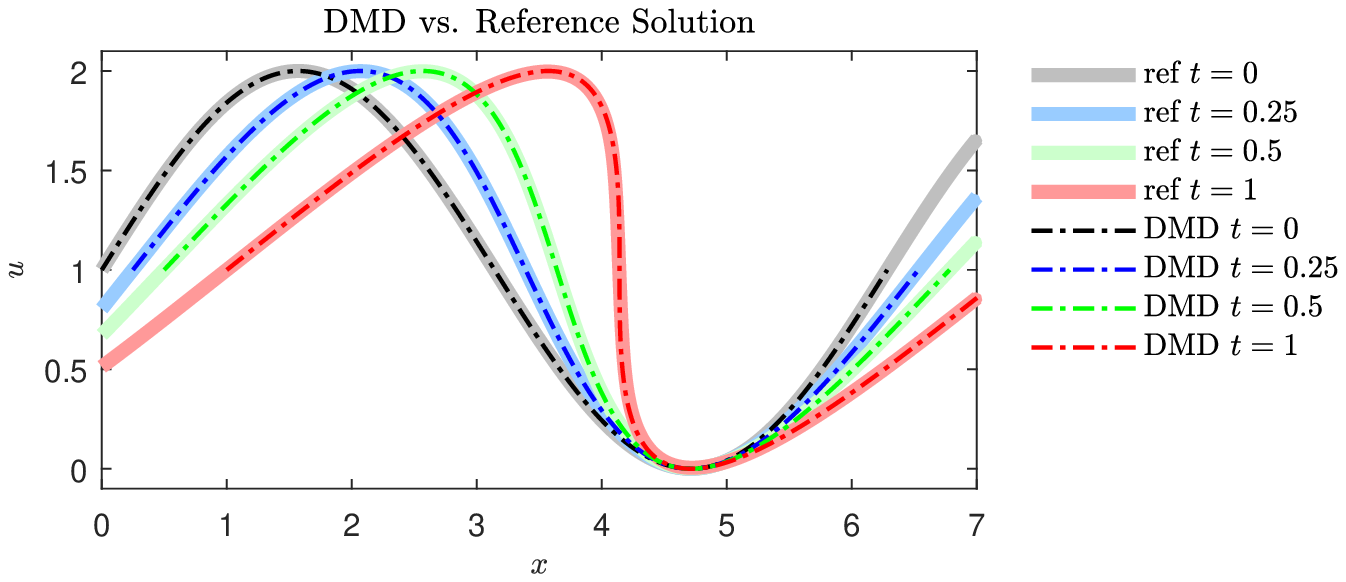}
\includegraphics{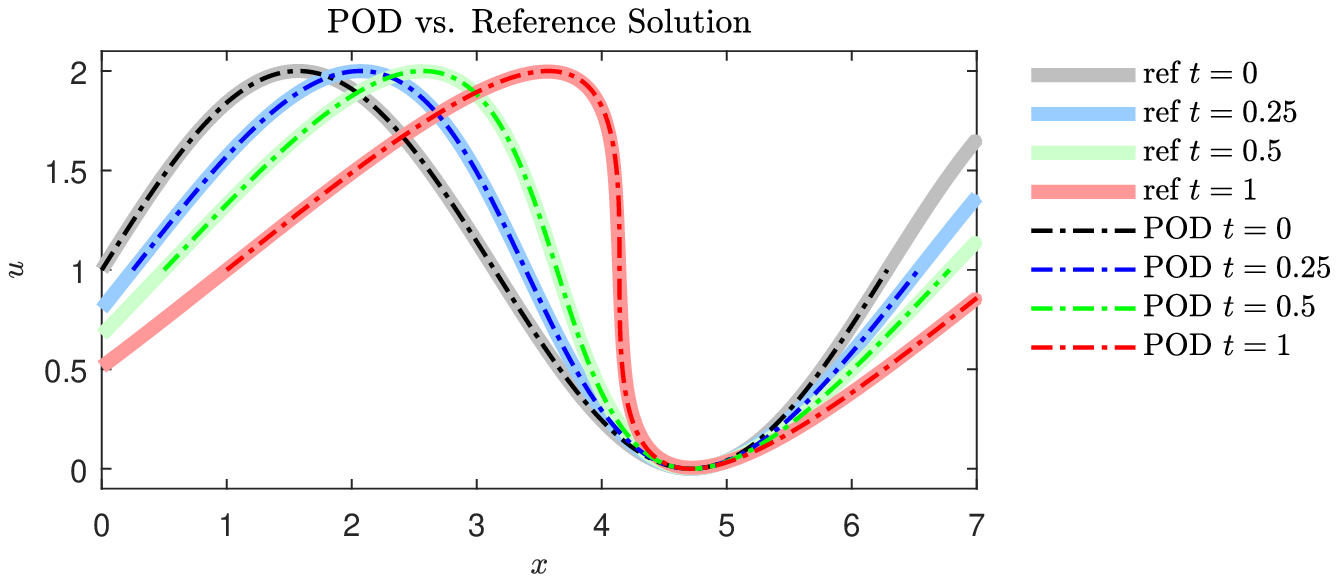}
\caption{Solutions of the inviscid Burgers equation, $u(x,t)$, alternatively obtained with the numerical method~\eqref{2-2} and the ROMs constructed via Lagrangian DMD (top) and Lagrangian POD (bottom). \label{fig:test3-vis} }
\end{figure}

\begin{figure}[htbp]
\centering
\includegraphics{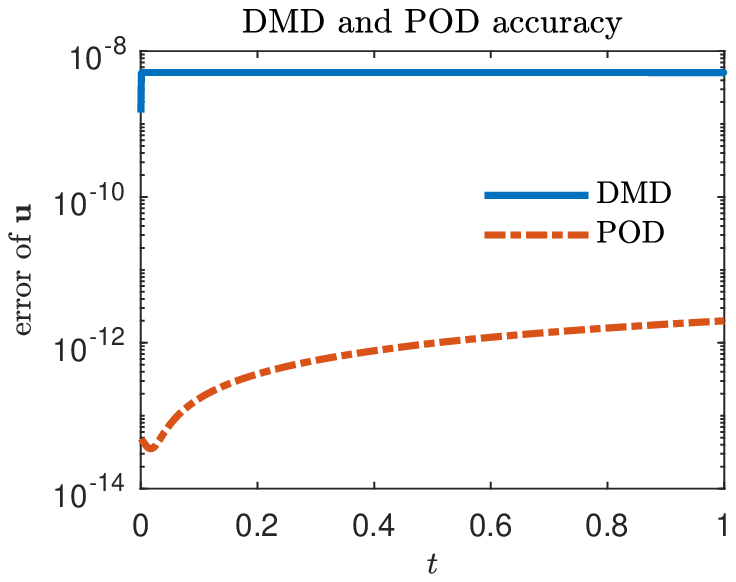}
\includegraphics{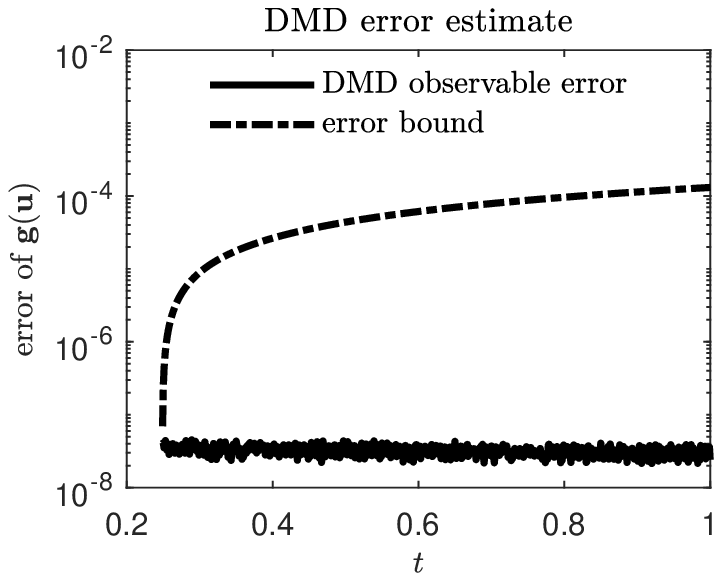}
\caption{Errors $\mathcal E$ of the Lagrangian DMD- and POD-based ROMs for the inviscid Burgers equation: error in reconstructing the state variable $u(x,t)$ (left) and its observables $\mathbf g(u)$ (right). The error bound for $\mathbf g(u)$ is derived in~\cite{lu2019error}. \label{fig:test3-err} }
\end{figure}

Figure~\ref{fig:test3-vis} provides a graphical illustration of the Lagrangian ROMs' ability to capture the system state dynamics in this nonlinear hyperbolic problem. Since the Lagrangian description treats first-order hyperbolic conservation laws, such as the inviscid Burgers equation, exactly via the method of characteristics, the addition of the particle trajectories $x(t)$ to the set of observables ensures that the Lagrangian POD and DMD are both of machine error accuracy (Figure~\ref{fig:test3-err}). Again the error bound serves as an indicator of accurate ROMs. 

The level-set method provides an alternative way to interpret the first-order hyperbolic conservation laws. In the Appendix, we report our experiments with the level-set DMD, which is essentially a Lagrangian DMD for two-dimensional linear advection equation.

\subsection{Viscous Burgers Equation}

The viscous Burgers equation is obtained from~\eqref{2-1} be setting $f \equiv u$ and $D \equiv 0.1$. Again, this equation is defined on $(x,t) \in (0,2\pi) \times (0,1]$ and is subject to the initial and boundary conditions~\eqref{4-3}.

\begin{figure}[htbp]
\centering
\includegraphics{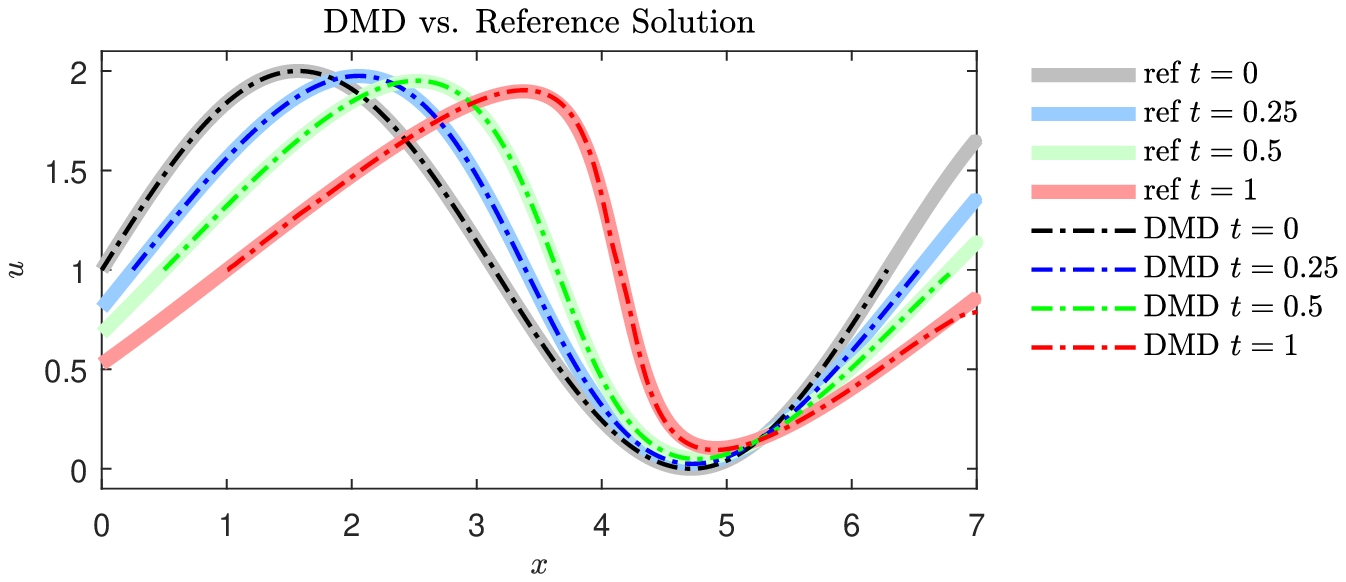}
\includegraphics{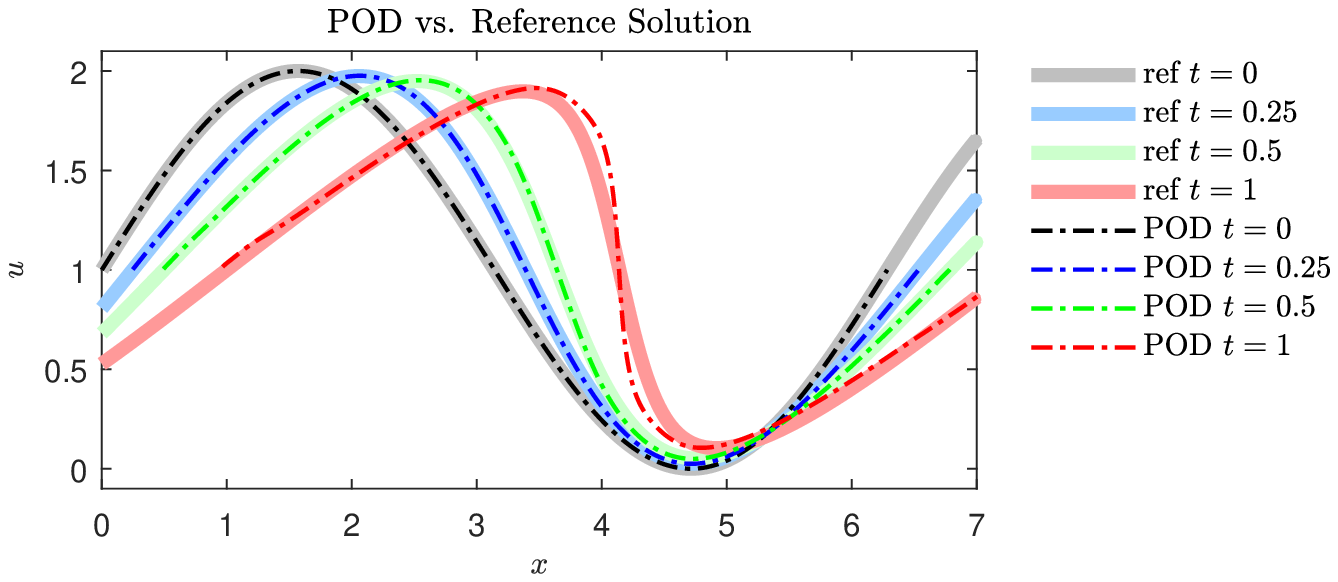}
\caption{Solutions of the viscous Burgers equation, $u(x,t)$, alternatively obtained with the numerical method~\eqref{2-2} and the ROMs constructed via Lagrangian DMD (top) and Lagrangian POD (bottom). \label{fig:test4-vis} }
\end{figure}

\begin{figure}[htbp]
\centering
\includegraphics{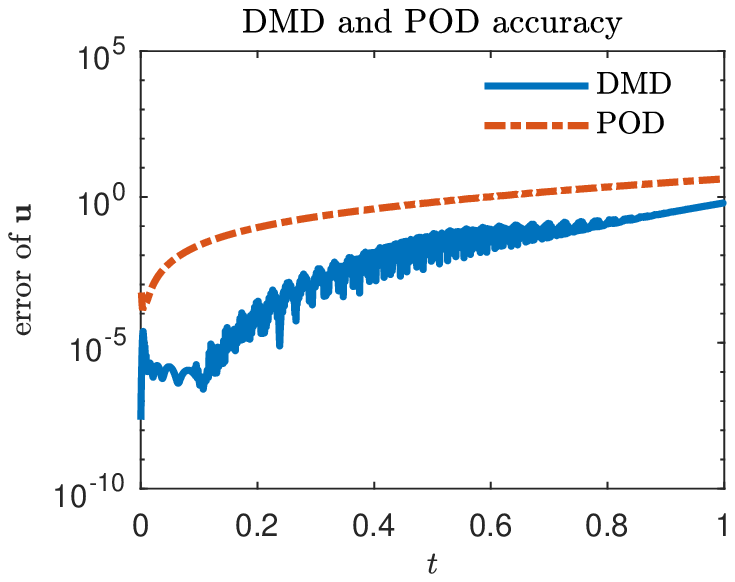}
\includegraphics{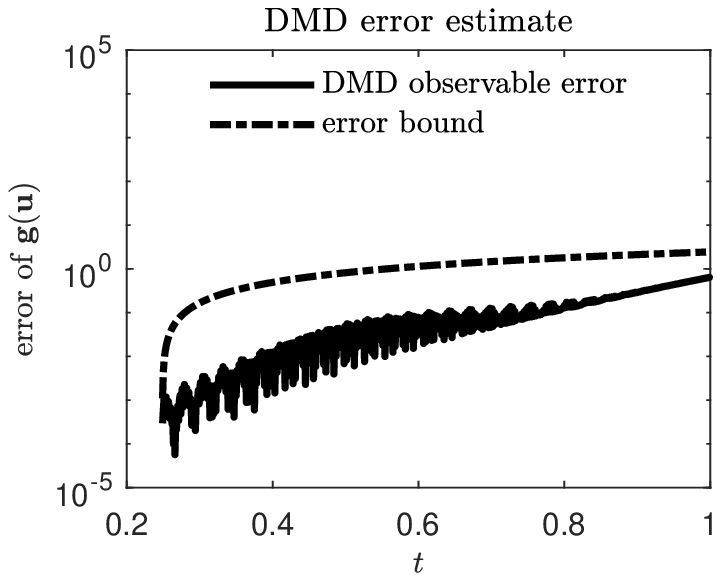}
\caption{Errors $\mathcal E$ of the Lagrangian DMD- and POD-based ROMs for the viscous Burgers equation: error in reconstructing the state variable $u(x,t)$ (left) and its observables $\mathbf g(u)$ (right). The error bound for $\mathbf g(u)$ is derived in~\cite{lu2019error}. \label{fig:test4-err} }
\end{figure}

For this nonlinear problem, Lagrangian DMD is visually more accurate that Lagrangian POD (Figure~\ref{fig:test4-vis}), especially at later times. This is confirmed by plotting the error $\mathcal E$ in Figure~\ref{fig:test4-err}. As mentioned in \cite{morton2018deep}, the Lagrangian grid might become distorted  (especially in the presence of large gradients $\partial_x u$) during the compressing process of ROM in the POD algorithm. The error estimate of the observable does a good job in evaluating the bound.

\subsection{Computational costs}

Table~\ref{tab:cost} collates the rank of the ROMs and the computational times of the HFM and the Lagrangian DMD- and POD-based ROMs.  In some cases, the SVD dominates the computational time of the ROM. Once the basis is constructed, the computation in the low rank subspace is much faster. This explains why the POD-based ROM of Test 3 takes more time to compute than the HFM. In other cases, the ROMs are much more efficient than the HFM computations. DMD is the most efficient methods because of its iteration-free nature.

\begin{table}[htbp]
\begin{center}
\begin{tabular}{ |c| c| c| c| c|}
\hline
 &Test 1& Test 2& Test 3& Test 4\\ 
 \hline
Rank truncation $r$ &3&10& 3& 14\\ 
 \hline
DMD computational time (sec) &0.114718& 0.046531& 0.048450&0.055641\\ 
 \hline
POD computational time (sec) &0.153869& 0.320905& 0.435342&0.795566 \\
  \hline
Eulerian HFM computational time (sec)&1.390459&29.782079&0.034519&55.713132\\
  \hline
  Lagrangian HFM computational time (sec)&0.023568 &27.020414&0.039063&55.246262\\
  \hline
\end{tabular}
\end{center}
\caption{Computational cost of the high-fidelity models and the corresponding Lagrangian DMD- and POD-based reduced-order models. Test 1 refers to advection problem; Test 2 to advection-diffusion problem; and Tests 3 and 4 to inviscid and viscous Burgers equations, respectively. \label{tab:cost} }
\end{table}

%%%%%%%%%%%%%%%%%%%%%%%%%%%%%%
\section{Conclusions}
\label{sec:concl}
%%%%%%%%%%%%%%%%%%%%%%%%%%%%%%

In this paper, we investigate the issue of translational problem for conventional proper orthogonal decomposition (POD) and dynamic mode decomposition (DMD) in the Eulerian framework. A new physic-aware DMD, based on the Lagrangian framework, is proposed to overcome the shortcomings of reduced order models (ROMs) of advection-dominated nonlinear phenomena. Characteristic lines, an important physical quantity in such systems, are taken into account in order to learn the Koopman operator of the underlying dynamics. The Lagrangian framework provides an optimal choice of observable functions for learning the Koopman operator. It allows one to construct a ROM in a relatively small subspace by using the DMD algorithm with satisfactory accuracy. Compared to the Lagrangian POD method, physics-aware DMD is more efficient computationally thanks to its iteration-free nature.

One possible direction for future work is to investigate the advection-diffusion system in Lagrangian coordinates~\cite{thiffeault2003advection}. Interpolation between Eulerian grid and Lagrangian grid will not be needed anymore but careful discretization of the diffusion operator will need to be handled. Existing numerical studies in Lagrangian coordinates and related methods~\cite{shashkov1999composite,hieber2005lagrangian} could be explored as guidelines of choosing physical observables in reduced order modeling.

All the numerical tests presented in this paper are shock-free. Once shock is formed, the Lagrangian formulation (\ref{3-1}) becomes invalid. Although one can still make the scheme (\ref{3-2}) work by numerical remedies, instability or unphysical solutions could appear when sharp gradients or shocks occur. The instability could become more severe in the compressed low-dimensional space~\cite{morton2018deep}. The modifications in \cite{morton2018deep} bypass this issue by compensating computational costs in projecting back to the Eulerian grid. From the perspective of physic-aware data-driven modeling, we realize that significant information like shock formation time, shock location and shock speed is not interpreted well enough from data. In another word, other quantities should be chosen as essential observables in order to learn the underlying Koopman operator. %In Part II of the physics-aware DMD study\cite{lu2019hodograph}, we are going to focus on shock informed DMD for ROM of conservation law.
We leave this line of research for a follow-up study.

\section*{Acknowledgements}
This work was supported in part by Defense Advanced Research Project Agency under award number 101513612 and by Air Force Office of Scientific Research under award number FA9550-17-1-0417.

%\clearpage

\begin{appendices}
\section{Level-set DMD for Hyperbolic Conservation Laws}

The level-set approach~\cite{tsai2003level} provides another way to interpret conservation laws. Supposed that a state variable $u(x,t)$ satisfies the one-dimensional conservation law  
\begin{equation}\label{cons}
\frac{\partial u}{\partial t} + f(u) \frac{\partial u}{\partial x} = 0
\end{equation}
with $f(u)\geq 0$. Its corresponding level-set formulation,
\begin{equation}\label{level}
\frac{\partial c}{\partial t} + f(y) \frac{\partial c}{\partial x} = 0,
\end{equation}
 is a linear two-dimensional transport equation for the dependent variable $c(x,y,t):\mathbb R^2\times \mathbb R^+\to \mathbb R$. Together with a Lipschitz-continuous initial function $c_0$, which embeds the initial data $u_0$ (see the example below), the zeroth-level set of $c(x,y,t)$, i.e., the $x-y$ contour of the solution to $c(x,y,t)=0$, gives the solution to the conservation law (\ref{cons}), $u(x,t)$. 
 
\begin{figure}[htbp]
\centering
\includegraphics{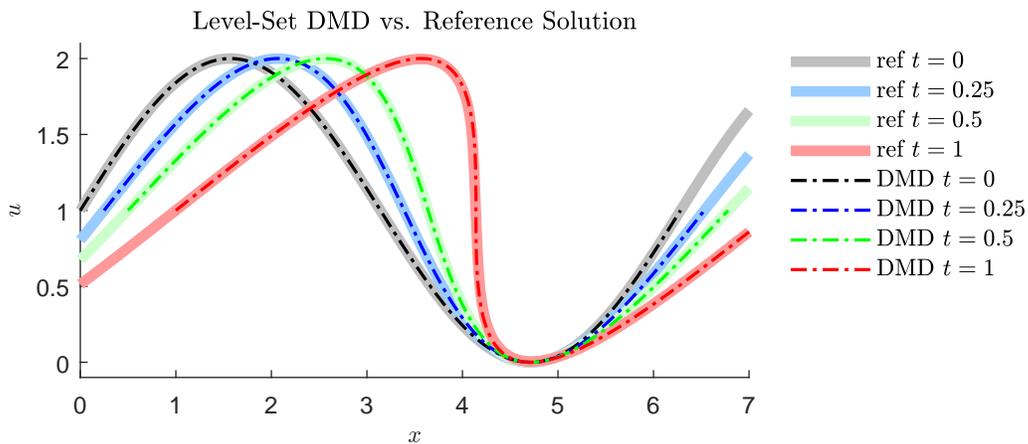}
\caption{Solutions of the inviscid Burgers equation, $u(x,t)$, alternatively obtained with the numerical method~\eqref{2-2} and the ROM constructed via the level-set DMD. \label{fig:level} }
\end{figure}

 By way of example, we consider the inviscid Burgers equation~\eqref{4-3}. Its level-set formulation is
\begin{equation}\label{level_burger}
\frac{\partial c}{\partial t} + y \frac{\partial c}{\partial x} = 0, \qquad c(x,y,0) = c_0(x,y) \equiv y - u_0.
\end{equation}
We apply Lagrangian DMD to construct a ROM for this two-dimensional linear advection equation from $m=250$ snapshots. Figure~\ref{fig:level} demonstrates that the ROM with the $r=3$ rank truncation approximates the HFM solution $u(x,t)$ with satisfactory accuracy. Although solving a two-dimensional linear problem takes more computational time than solving the nonlinear one-dimensional problem in this case, the level-set DMD provides another venue for investigation of physics-aware DMD that might have efficient applications to other problems.

\end{appendices}

\renewcommand\refname{Reference}
\bibliography{LagrangianDMD}

\end{document}